\documentclass[12pt]{article}
\usepackage{amsfonts,amssymb,amsmath}

\font\smallit=cmti10

\newtheorem{theorem}{Theorem}
\newtheorem{corollary}[theorem]{Corollary}
\newtheorem{lemma}[theorem]{Lemma}
\newtheorem{proposition}[theorem]{Proposition}
\newtheorem{remark}[theorem]{Remark}

\begin{document}

\begin{center}
{\bf A NOTE ON BOOLEAN LATTICES AND FAREY SEQUENCES} \vskip 20pt {\bf Andrey O. Matveev}\\ {\smallit
Data-Center Co., RU-620034, Ekaterinburg, P.O.~Box~5, Russian~Federation}\\ {\tt aomatveev@dc.ru,
aomatveev@hotmail.com}
\end{center}
\vskip 30pt

\centerline{\bf Abstract}

\noindent We establish monotone bijections between the Farey sequences of order~$m$ and the halfsequences of
Farey subsequences associated with the rank~$m$ elements of the Boolean lattice of subsets of a~$2m$-set. We
also present a few related combinatorial identities.

\noindent {\small Subject class: 05A19, 11B65.}


\thispagestyle{empty} \baselineskip=15pt \vskip 30pt

\section*{\normalsize 1. Introduction}

The {\em Farey sequence of order $n$}, denoted by $\mathcal{F}_n$, is the ascending sequence of irreducible
fractions $\tfrac{h}{k}\in\mathbb{Q}$ with $0\leq h\leq k\leq n$, see, e.g.,~[2, Chapter~27], [3, \S{}3], [4,
Chapter~4], \mbox{[5, Chapter~III],} [8, Chapter~6], [9, Chapter~5]; their numerators and denominators are
presented in sequences A006842 and A006843 in Sloane's {\sl On-Line Encyclopedia of Integer Sequences}. For
example,
\begin{equation*}
\mathcal{F}_6= \left(
\tfrac{0}{1}<\tfrac{1}{6}<\tfrac{1}{5}<\tfrac{1}{4}<\tfrac{1}{3}<\tfrac{2}{5}<
\tfrac{1}{2}<\tfrac{3}{5}<\tfrac{2}{3}<\tfrac{3}{4}<\tfrac{4}{5}<\tfrac{5}{6}<
\tfrac{1}{1} \right)\ .
\end{equation*}

For any integer~$m$, $0<m<n$, the ascending sets
\begin{equation}
\label{eq:10} \left(\tfrac{h}{k}\in\mathcal{F}_n:\ h\leq m\right)
\end{equation}
are interesting Farey subsequences~[1].

Let $C$ be a finite set of cardinality $n:=|C|$ greater than or
equal to two, and $A$ its proper subset; $m:=|A|$. Denote the
Boolean lattice of subsets of~$C$ by $\mathbb{B}(n)$; the empty
set is denoted by $\hat{0}$, and the family of $l$-element subsets
of $C$ is denoted  by $\mathbb{B}(n)^{(l)}$. Let
$\gcd(\cdot,\cdot)$ denote the greatest common divisor of two
integers. The ascending sequence of fractions
\begin{equation*}
\begin{split}
\mathcal{F}\bigl(\mathbb{B}(n),m\bigr):&=\left(\tfrac{|B\cap
A|}{\gcd(|B\cap A|,|B|)}\!\!\Bigm/\!\!\!\tfrac{|B|}{\gcd(|B\cap
A|,|B|)}:\ B\subseteq C,\
|B|>0\right)\\&=\left(\tfrac{h}{k}\in\mathcal{F}_n:\ h\leq m,\
k-h\leq n-m\right)\ ,
\end{split}
\end{equation*}
considered in~[7], has the properties very similar to those of the standard Farey sequence $\mathcal{F}_n$
and of Farey subsequence~(\ref{eq:10}).

The Farey subsequences
\begin{equation*}
\mathcal{F}\bigl(\mathbb{B}(2m),m\bigr):=\left(\tfrac{h}{k}\in\mathcal{F}_{2m}:\
h\leq m,\ k-h\leq m\right)
\end{equation*}
arise in analysis of decision-making problems~[6]. One of such subsequences is
\begin{multline*}
\mathcal{F}\bigl(\mathbb{B}(12),6\bigr)=
\bigl(\tfrac{0}{1}<\tfrac{1}{7}<\tfrac{1}{6}<\tfrac{1}{5}<\tfrac{1}{4}<\tfrac{2}{7}<\tfrac{1}{3}
<\tfrac{3}{8}<\tfrac{2}{5}<\tfrac{3}{7}<\tfrac{4}{9}<\tfrac{5}{11}<\tfrac{1}{2}\\
<\tfrac{6}{11}<\tfrac{5}{9}<\tfrac{4}{7}<\tfrac{3}{5}<\tfrac{5}{8}<\tfrac{2}{3}
<\tfrac{5}{7}<\tfrac{3}{4}<\tfrac{4}{5}<\tfrac{5}{6}<\tfrac{6}{7}<\tfrac{1}{1}\bigr)\
.
\end{multline*}

The fractions in the above-mentioned Farey (sub)sequences are indexed starting with zero.

In Theorem~\ref{prop:6} of this note we establish the connection between the standard Farey sequence
$\mathcal{F}_m$ and the halfsequences of $\mathcal{F}\bigl(\mathbb{B}(2m),m\bigr)$.

\vskip 30pt

\section*{\normalsize 2. The Farey Subsequence $\mathcal{F}\bigl(\mathbb{B}(n),m\bigr)$}

Recall that the map $\mathcal{F}_n\to\mathcal{F}_n$, which sends a
fraction $\tfrac{h}{k}$ to $\tfrac{k-h}{k}$, is order-reversing
and bijective. The sequences
$\mathcal{F}\bigl(\mathbb{B}(n),m\bigr)$ and
$\mathcal{F}\bigl(\mathbb{B}(n),n-m\bigr)$ have an analogous
property:

\begin{lemma}{\rm[7]}
\label{prop:3} The map
\begin{equation}
\label{eq:9}
\mathcal{F}\bigl(\mathbb{B}(n),m\bigr)\to\mathcal{F}\bigl(\mathbb{B}(n),n-m\bigr)\
,\ \ \ \tfrac{h}{k}\mapsto\tfrac{k-h}{k}\ ,
\end{equation}
is order-reversing and bijective.
\end{lemma}

If we write the fractions $\tfrac{h}{k}\in\mathbb{Q}$ as the
column vectors $\left[\begin{smallmatrix}\!h\!\\
\!k\!\end{smallmatrix}\right]\in\mathbb{Z}^2$, then
map~(\ref{eq:9}) can be thought of as the map
\begin{equation*}
\left[\begin{smallmatrix}\!h\!\\
\!k\!\end{smallmatrix}\right]\mapsto\left[\begin{smallmatrix}-1&1\\0&1\end{smallmatrix}\right]\cdot
\left[\begin{smallmatrix}\!h\!\\ \!k\!\end{smallmatrix}\right]\ .
\end{equation*}

Let $a'\in\mathbb{B}(n)$ and $0<m:=\rho(a')<n$, where $\rho(a')$
denotes the poset rank of~$a'$ in $\mathbb{B}(n)$. For a subset
$A\subset\mathbb{B}(n)$, let $\mathfrak{I}(A)$ and
$\mathfrak{F}(A)$ denote the order ideal and filter in
$\mathbb{B}(n)$, generated by $A$, respectively. The subposet
\mbox{$\mathfrak{F}\bigl(\mathfrak{I}(a')\cap\mathbb{B}(n)^{(1)}\bigr)$,}
of cardinality \mbox{$2^n-2^{n-m}$,} can be partitioned in the
following way:
\begin{multline*}
\mathfrak{F}\bigl(\mathfrak{I}(a')\cap\mathbb{B}(n)^{(1)}\bigr)= \bigl(\mathfrak{I}(a')-\{\hat{0}\}\bigr)\ \
\dot\cup\\ \dot{\bigcup_{\substack{\frac{h}{k}\in\mathcal{F}(\mathbb{B}(n),m):\\
\frac{0}{1}<\frac{h}{k}<\frac{1}{1}}}}\ \ \dot{\bigcup_{1\leq s\leq \left\lfloor\min\left\{\frac{m}{h},\
\frac{n-m}{k-h}\right\}\right\rfloor}} \Bigl(\mathbb{B}(n)^{(s\cdot k)}\cap
\bigl(\mathfrak{F}\bigl(\mathfrak{I}(a')\cap\mathbb{B}(n)^{(s\cdot h)}\bigr)-
\mathfrak{F}\bigl(\mathfrak{I}(a')\cap\mathbb{B}(n)^{(s\cdot h +1)}\bigr)\bigr)\Bigr)\ ,
\end{multline*}
where
\begin{equation*}
\bigl|\mathbb{B}(n)^{(s\cdot k)}\cap
\bigl(\mathfrak{F}\bigl(\mathfrak{I}(a')\cap\mathbb{B}(n)^{(s\cdot
h)}\bigr)-
\mathfrak{F}\bigl(\mathfrak{I}(a')\cap\mathbb{B}(n)^{(s\cdot
h+1)}\bigr)\bigr)\bigr|= \tbinom{m}{s\cdot
h}\tbinom{n-m}{s\cdot(k-h)}\ .
\end{equation*}
Since $|\mathfrak{I}(a')-\{\hat{0}\}|=2^m-1$, we obtain
\begin{equation*}
2^n-2^{n-m}=2^m-1
+\sum_{\substack{\frac{h}{k}\in\mathcal{F}(\mathbb{B}(n),m):\\
\frac{0}{1}<\frac{h}{k}<\frac{1}{1}}}\ \ \sum_{1\leq s\leq
\left\lfloor\min\left\{\frac{m}{h},\
\frac{n-m}{k-h}\right\}\right\rfloor}\tbinom{m}{s\cdot
h}\tbinom{n-m}{s\cdot(k-h)}\ .
\end{equation*}

If $a''\in\mathbb{B}(n)$ and $\rho(a'')=n-m$, then
Lemma~\ref{prop:3} implies
\begin{equation*}
2^n-2^m=2^{n-m}-1
+\sum_{\substack{\frac{h}{k}\in\mathcal{F}(\mathbb{B}(n),n-m):\\
\frac{0}{1}<\frac{h}{k}<\frac{1}{1}}}\ \ \sum_{1\leq s\leq
\left\lfloor\min\left\{\frac{n-m}{h},\
\frac{m}{k-h}\right\}\right\rfloor}\tbinom{n-m}{s\cdot
h}\tbinom{m}{s\cdot(k-h)}\ ,
\end{equation*}
and we come to the following conclusion:

\begin{proposition}  \label{prop:4} Fractions from the Farey subsequences
$\mathcal{F}(\mathbb{B}(n),m)$ and
$\mathcal{F}(\mathbb{B}(n),n-m)$ satisfy the equality:
\begin{align*}
\quad&\phantom{=}\,\,
\sum_{\substack{\frac{h}{k}\in\mathcal{F}(\mathbb{B}(n),m):\\
\frac{0}{1}<\frac{h}{k}<\frac{1}{1}}}\ \sum_{1\leq s\leq
\left\lfloor\min\left\{\frac{m}{h},\
\frac{n-m}{k-h}\right\}\right\rfloor}\tbinom{m}{s\cdot
h}\tbinom{n-m}{s\cdot(k-h)}\\&=
\sum_{\substack{\frac{h}{k}\in\mathcal{F}(\mathbb{B}(n),n-m):\\
\frac{0}{1}<\frac{h}{k}<\frac{1}{1}}}\ \ \sum_{1\leq s\leq
\left\lfloor\min\left\{\frac{n-m}{h},\
\frac{m}{k-h}\right\}\right\rfloor}\tbinom{n-m}{s\cdot
h}\tbinom{m}{s\cdot(k-h)}\\&= 2^n-2^m-2^{n-m}+1\ .
\end{align*}
\end{proposition}

\vskip 30pt

\section*{\normalsize 3. The Farey Subsequence $\mathcal{F}\bigl(\mathbb{B}(2m),m\bigr)$}

Denote the left and right halfsequences of
$\mathcal{F}\bigl(\mathbb{B}(2m),m\bigr)$ by
\begin{align*}
\mathcal{F}^{\leq\frac{1}{2}}\!\bigl(\mathbb{B}(2m),m\bigr):&=
\left(f\in\mathcal{F}\bigl(\mathbb{B}(2m),m\bigr):\
f\leq\tfrac{1}{2}\right) \intertext{and}
\mathcal{F}^{\geq\frac{1}{2}}\!\bigl(\mathbb{B}(2m),m\bigr):&=
\left(f\in\mathcal{F}\bigl(\mathbb{B}(2m),m\bigr):\
f\geq\tfrac{1}{2}\right)\ ,
\end{align*}
respectively.
\begin{lemma}{\rm[6]}
\label{prop:1} The maps
\begin{align*}
\mathcal{F}\bigl(\mathbb{B}(2m),m\bigr)&\to\mathcal{F}\bigl(\mathbb{B}(2m),m\bigr)\
, &  \tfrac{h}{k}&\mapsto\tfrac{k-h}{k}\ , &
\left[\begin{smallmatrix}\!h\!\\
\!k\!\end{smallmatrix}\right]&\mapsto\left[\begin{smallmatrix}-1&1\\0&1\end{smallmatrix}\right]\cdot
\left[\begin{smallmatrix}\!h\!\\ \!k\!\end{smallmatrix}\right]\
,\\
\mathcal{F}^{\leq\frac{1}{2}}\!\bigl(\mathbb{B}(2m),m\bigr)&\to
\mathcal{F}^{\leq\frac{1}{2}}\!\bigl(\mathbb{B}(2m),m\bigr)\ , &
\tfrac{h}{k}&\mapsto\tfrac{k-2h}{2k-3h}\ , &
\left[\begin{smallmatrix}\!h\!\\
\!k\!\end{smallmatrix}\right]&\mapsto\left[\begin{smallmatrix}-2&1\\-3&2\end{smallmatrix}\right]\cdot
\left[\begin{smallmatrix}\!h\!\\ \!k\!\end{smallmatrix}\right]\
,\\ \intertext{and}
\mathcal{F}^{\geq\frac{1}{2}}\!\bigl(\mathbb{B}(2m),m\bigr)&\to
\mathcal{F}^{\geq\frac{1}{2}}\!\bigl(\mathbb{B}(2m),m\bigr)\ , &
\tfrac{h}{k}&\mapsto\tfrac{h}{3h-k}\ , &
\left[\begin{smallmatrix}\!h\!\\
\!k\!\end{smallmatrix}\right]&\mapsto\left[\begin{smallmatrix}1&0\\3&-1\end{smallmatrix}\right]\cdot
\left[\begin{smallmatrix}\!h\!\\ \!k\!\end{smallmatrix}\right]\ ,
\end{align*}
are order-reversing and bijective.
\end{lemma}

\begin{corollary}
\label{prop:2} The maps
\begin{align*}
\mathcal{F}^{\leq\frac{1}{2}}\!\bigl(\mathbb{B}(2m),m\bigr)&\to
\mathcal{F}^{\geq\frac{1}{2}}\!\bigl(\mathbb{B}(2m),m\bigr)\ , &
\tfrac{h}{k}&\mapsto\tfrac{k-h}{2k-3h}\ , &
\left[\begin{smallmatrix}\!h\!\\
\!k\!\end{smallmatrix}\right]&\mapsto\left[\begin{smallmatrix}-1&1\\-3&2\end{smallmatrix}\right]\cdot
\left[\begin{smallmatrix}\!h\!\\ \!k\!\end{smallmatrix}\right]\
,\\ \intertext{and}
\mathcal{F}^{\geq\frac{1}{2}}\!\bigl(\mathbb{B}(2m),m\bigr)&\to
\mathcal{F}^{\leq\frac{1}{2}}\!\bigl(\mathbb{B}(2m),m\bigr)\ , &
\tfrac{h}{k}&\mapsto\tfrac{2h-k}{3h-k}\ , &
\left[\begin{smallmatrix}\!h\!\\
\!k\!\end{smallmatrix}\right]&\mapsto\left[\begin{smallmatrix}2&-1\\3&-1\end{smallmatrix}\right]\cdot
\left[\begin{smallmatrix}\!h\!\\ \!k\!\end{smallmatrix}\right]\ ,
\end{align*}
are order-preserving and bijective.
\end{corollary}

Let $f_{t^1_3}$, $f_{t^1_2}$, $f_{t^2_3}$,
$f_{t^1_1}\in\mathcal{F}\bigl(\mathbb{B}(2m),m\bigr)$, $m>1$,
where
\begin{equation*}
f_{t^1_3}:=\tfrac{1}{3}\ ,\ \ \ f_{t^1_2}:=\tfrac{1}{2}\ ,\ \ \
f_{t^2_3}:=\tfrac{2}{3}\ ,\ \ \ f_{t^1_1}:=\tfrac{1}{1}\ ;
\end{equation*}
then Lemma~{\rm\ref{prop:1}} and Corollary~{\rm\ref{prop:2}} imply
\begin{equation*}
t^1_2=2t^1_3\ ,\ \ \ t^2_3=3t^1_3\ ,\ \ \ t^1_1=4t^1_3\ .
\end{equation*}
This in particular means that the number
$\left|\mathcal{F}\bigl(\mathbb{B}(2m),m\bigr)\right|-1=:t^1_1$ is
divisible by four.

\vskip 30pt

\section*{\normalsize 4. The Farey Sequence $\mathcal{F}_m$ and the Farey Subsequence
$\mathcal{F}\bigl(\mathbb{B}(2m),m\bigr)$}

Let $h$ be a positive integer, and $[i,l]:=\{j:\ i\leq j\leq l\}$
an interval of positive integers. Let
\begin{equation*}
\phi(h;[i,l]):=\bigl|\bigl\{j\in[i,l]:\ \gcd(h,j)=1\bigr\}\bigr|\
;
\end{equation*}
thus, $\phi(h;[1,h])$ is the {\em Euler $\phi$-function}. Recall
that for a nonempty interval of positive integers $[i'+1,i'']$ it
holds $\phi(h;[i'+1,i''])=\sum_{d\in[1,\min\{i'',h\}]:\ d|h}\,
\overline{\mu}(d)$
$\cdot\left(\left\lfloor\tfrac{i''}{d}\right\rfloor
-\left\lfloor\tfrac{i'}{d}\right\rfloor\right)$, where $d|h$ means
that $d$ divides $h$, and $\overline{\mu}(\cdot)$ stands for the
{\em M\"{o}bius function}: $\overline{\mu}(1):=1$; if $p^2|d$, for
some prime $p$, then $\overline{\mu}(d):=0$; if $d=p_1 p_2\cdots
p_s$ is the product of distinct primes $p_1,p_2,\ldots,p_s$, then
$\overline{\mu}(d):=(-1)^s$.

Let $m$ be an integer, $m>1$. For every integer $h$, $1\leq h\leq
m$, we have
\begin{multline*}
\bigr|\bigl\{\tfrac{h}{k}\in\mathcal{F}\bigl(\mathbb{B}(2m),m\bigr):\
\tfrac{h}{k}<\tfrac{1}{2}\bigr\}\bigr|=\phi(h;[2h+1,h+m])\\=\sum_{d\in[1,h]:\
d|h}\overline{\mu}(d)\cdot\left(\left\lfloor\tfrac{h+m}{d}\right\rfloor
-\tfrac{2h}{d}\right)=\sum_{d\in[1,h]:\
d|h}\overline{\mu}(d)\cdot\left(\left\lfloor\tfrac{m}{d}\right\rfloor
-\tfrac{h}{d}\right)\\=\phi(h;[h+1,m])=\bigl|\bigl\{\tfrac{h}{k}\in\mathcal{F}_m:\
\tfrac{h}{k}<\tfrac{1}{1}\bigr\}\bigr|\ ;
\end{multline*}
hence, the sequences
$\mathcal{F}^{\leq\frac{1}{2}}\!\bigl(\mathbb{B}(2m),m\bigr)$ and
$\mathcal{F}_m$ are of the same cardinality. Noticing that
fractions $\tfrac{h_j}{k_j}$ and $\tfrac{h_{j+1}}{k_{j+1}}$ are
consecutive in $\mathcal{F}_m$ if and only if the fractions
$\tfrac{h_j}{k_j+h_j}$ and $\tfrac{h_{j+1}}{k_{j+1}+h_{j+1}}$ are
consecutive in
$\mathcal{F}^{\leq\frac{1}{2}}\!\bigl(\mathbb{B}(2m),m\bigr)$, we
arrive, with the help of Lemma~{\rm\ref{prop:1}} and
Corollary~{\rm\ref{prop:2}}, at the following conclusion:

\begin{theorem}\label{prop:6} Let $m$ be an integer, $m>1$.
The maps
\begin{align*}
\mathcal{F}^{\leq\frac{1}{2}}\!\bigl(\mathbb{B}(2m),m\bigr)&\to\mathcal{F}_m\
, & \tfrac{h}{k}&\mapsto\tfrac{h}{k-h}\ , &
\left[\begin{smallmatrix}\!h\!\\
\!k\!\end{smallmatrix}\right]&\mapsto\left[\begin{smallmatrix}1&0\\-1&1\end{smallmatrix}\right]\cdot
\left[\begin{smallmatrix}\!h\!\\ \!k\!\end{smallmatrix}\right]\
,\intertext{and}
\mathcal{F}_m&\to\mathcal{F}^{\leq\frac{1}{2}}\!\bigl(\mathbb{B}(2m),m\bigr)\
, & \tfrac{h}{k}&\mapsto\tfrac{h}{k+h}\ , &
\left[\begin{smallmatrix}\!h\!\\
\!k\!\end{smallmatrix}\right]&\mapsto\left[\begin{smallmatrix}1&0\\1&1\end{smallmatrix}\right]\cdot
\left[\begin{smallmatrix}\!h\!\\ \!k\!\end{smallmatrix}\right]\ ,
\end{align*}
are order-preserving and bijective.

The maps
\begin{align*}
\mathcal{F}^{\geq\frac{1}{2}}\!\bigl(\mathbb{B}(2m),m\bigr)&\to\mathcal{F}_m\
, & \tfrac{h}{k}&\mapsto\tfrac{k-h}{h}\ , &
\left[\begin{smallmatrix}\!h\!\\
\!k\!\end{smallmatrix}\right]&\mapsto\left[\begin{smallmatrix}-1&1\\1&0\end{smallmatrix}\right]\cdot
\left[\begin{smallmatrix}\!h\!\\ \!k\!\end{smallmatrix}\right]\
,\intertext{and}
\mathcal{F}_m&\to\mathcal{F}^{\geq\frac{1}{2}}\!\bigl(\mathbb{B}(2m),m\bigr)\
, & \tfrac{h}{k}&\mapsto\tfrac{k}{k+h}\ , &
\left[\begin{smallmatrix}\!h\!\\
\!k\!\end{smallmatrix}\right]&\mapsto\left[\begin{smallmatrix}0&1\\1&1\end{smallmatrix}\right]\cdot
\left[\begin{smallmatrix}\!h\!\\ \!k\!\end{smallmatrix}\right]\ ,
\end{align*}
are order-reversing and bijective.
\end{theorem}

Direct counting gives $\left|\mathcal{F}\bigl(\mathbb{B}(2),1\bigr)\right|=3$ and
$\left|\mathcal{F}\bigl(\mathbb{B}(4),2\bigr)\right|=5$.

Since $|\mathcal{F}_m|-1=\tfrac{1}{2}\sum_{d\geq
1}\overline{\mu}(d)\cdot\left\lfloor\tfrac{m}{d}\right\rfloor\cdot \left\lfloor\tfrac{m}{d}+1\right\rfloor$
(see, e.g.,~[4, \S{}4.9]), Theorem~\ref{prop:6} implies that for $m>1$ we have
\begin{equation*}
\left|\mathcal{F}\bigl(\mathbb{B}(2m),m\bigr)\right|-1=\sum_{d\geq
1}\overline{\mu}(d)\cdot\left\lfloor\tfrac{m}{d}\right\rfloor\cdot
\left\lfloor\tfrac{m}{d}+1\right\rfloor\ .
\end{equation*}

By means of Theorem~\ref{prop:6}, the descriptions of sequences $\mathcal{F}_m$ and
$\mathcal{F}\bigl(\mathbb{B}(2m),m\bigr)$ supplement each other. For example, consider a fraction
$\tfrac{h}{k}\in\mathcal{F}_m-\left\{\tfrac{0}{1},\tfrac{1}{1}\right\}$. If $x_0$ is the integer such that
$hx_0\equiv -1\pmod{k}$ and $m-k+1\leq x_0\leq m$, then it is known (see, e.g.,~[2, \S{}27.1]) that the
fraction $\tfrac{hx_0+1}{k}\!\!\!\Bigm/\!\!\! x_0$ succeeds the fraction $\tfrac{h}{k}$ in $\mathcal{F}_m$.
Similarly, if $x_0$ is the integer such that $hx_0\equiv 1\pmod{k}$ and $m-k+1\leq x_0\leq m$, then the
fraction $\tfrac{hx_0-1}{k}\!\!\!\Bigm/\!\!\! x_0$ precedes $\tfrac{h}{k}$ in $\mathcal{F}_m$.
Theorem~\ref{prop:6} leads to an analogous statement:

\newpage

\begin{remark} Let $m$ be an integer, $m>1$.
\begin{itemize}
\item[\rm(i)]
Let $\tfrac{h}{k}\in\mathcal{F}\bigl(\mathbb{B}(2m),m\bigr)$.
Suppose that $\tfrac{0}{1}<\tfrac{h}{k}\leq\tfrac{1}{2}$. Let
$x_0$ be the integer such that $hx_0\equiv 1\pmod{(k-h)}$ and
$m-k+h+1\leq x_0\leq m$. The fraction
\begin{equation*}
\tfrac{hx_0-1}{k-h}\!\!\Bigm/\!\!\tfrac{kx_0-1}{k-h}
\end{equation*}
precedes $\tfrac{h}{k}$ in
$\mathcal{F}\bigl(\mathbb{B}(2m),m\bigr)$.
\item[\rm(ii)]
Let $\tfrac{h}{k}\in\mathcal{F}\bigl(\mathbb{B}(2m),m\bigr)$.
Suppose that $\tfrac{0}{1}\leq\tfrac{h}{k}<\tfrac{1}{2}$. Let
$x_0$ be the integer such that $hx_0\equiv -1\pmod{(k-h)}$ and
$m-k+h+1\leq x_0\leq m$. The fraction
\begin{equation*}
\tfrac{hx_0+1}{k-h}\!\!\Bigm/\!\!\tfrac{kx_0+1}{k-h}
\end{equation*}
succeeds $\tfrac{h}{k}$ in
$\mathcal{F}\bigl(\mathbb{B}(2m),m\bigr)$.
\end{itemize}
\end{remark}

Proposition~\ref{prop:4} can be reformulated in the case where $n:=2m$, with the help of the bijections
mentioned in Lemma~\ref{prop:1} and Corollary~\ref{prop:2}, in several ways which we now summarize:

\begin{proposition}\label{prop:12}
Let $m$ be an integer, $m>1$. The following combinatorial
identities hold for fractions from the Farey subsequence
$\mathcal{F}\bigl(\mathbb{B}(2m),m\bigr)$:
\begin{itemize}
\item[\rm(i)]
\begin{equation*}
\sum_{\substack{\frac{h}{k}\in\mathcal{F}(\mathbb{B}(2m),m):\\
\frac{0}{1}<\frac{h}{k}<\frac{1}{1}}}\ \sum_{1\leq s\leq
\left\lfloor\min\left\{\frac{m}{h},\
\frac{m}{k-h}\right\}\right\rfloor}\tbinom{m}{s\cdot
h}\tbinom{m}{s\cdot(k-h)}=2^{2m}-2^{m+1}+1\ .
\end{equation*}
\item[\rm(ii)]
\begin{equation*}
\begin{split}
\sum_{\substack{\frac{h}{k}\in\mathcal{F}(\mathbb{B}(2m),m):\\
\frac{0}{1}<\frac{h}{k}<\frac{1}{2}}}\ \sum_{1\leq s\leq
\left\lfloor\frac{m}{k-h}\right\rfloor}\tbinom{m}{s\cdot
h}\tbinom{m}{s\cdot(k-h)}&=
\sum_{\substack{\frac{h}{k}\in\mathcal{F}(\mathbb{B}(2m),m):\\
\frac{1}{2}<\frac{h}{k}<\frac{1}{1}}}\ \sum_{1\leq s\leq
\left\lfloor\frac{m}{h}\right\rfloor}\tbinom{m}{s\cdot
h}\tbinom{m}{s\cdot(k-h)}\\
&=2^{2m-1}-2^m-\tfrac{1}{2}\tbinom{2m}{m}+1\ .
\end{split}
\end{equation*}
\item[\rm(iii)]
\begin{equation*}
\begin{split}
\quad&\phantom{=}\sum_{\substack{\frac{h}{k}\in\mathcal{F}(\mathbb{B}(2m),m):\\
\frac{0}{1}<\frac{h}{k}<\frac{1}{3}}}\ \sum_{1\leq s\leq
\left\lfloor\frac{m}{k-h}\right\rfloor}\tbinom{m}{s\cdot(k-h)}\left(\tbinom{m}{s\cdot
h}+\tbinom{m}{s\cdot (k-2h)}\right)\\&=
\sum_{\substack{\frac{h}{k}\in\mathcal{F}(\mathbb{B}(2m),m):\\
\frac{1}{3}<\frac{h}{k}<\frac{1}{2}}}\ \sum_{1\leq s\leq
\left\lfloor\frac{m}{k-h}\right\rfloor}\tbinom{m}{s\cdot(k-h)}\left(\tbinom{m}{s\cdot
h}+\tbinom{m}{s\cdot (k-2h)}\right)
\end{split}
\end{equation*}
\begin{equation*}
\begin{split}
&=\sum_{\substack{\frac{h}{k}\in\mathcal{F}(\mathbb{B}(2m),m):\\ \frac{1}{2}<\frac{h}{k}<\frac{2}{3}}}\
\sum_{1\leq s\leq \left\lfloor\frac{m}{h}\right\rfloor}\tbinom{m}{s\cdot
h}\left(\tbinom{m}{s\cdot(k-h)}+\tbinom{m}{s\cdot(2h-k)}\right)\\&=
\sum_{\substack{\frac{h}{k}\in\mathcal{F}(\mathbb{B}(2m),m):\\ \frac{2}{3}<\frac{h}{k}<\frac{1}{1}}}\
\sum_{1\leq s\leq \left\lfloor\frac{m}{h}\right\rfloor}\tbinom{m}{s\cdot
h}\left(\tbinom{m}{s\cdot(k-h)}+\tbinom{m}{s\cdot(2h-k)}\right)\\ &=2^{2m-1}-2^m-\tfrac{1}{2}\tbinom{2m}{m}-
\sum_{1\leq t\leq \left\lfloor\frac{m}{2}\right\rfloor}\tbinom{m}{2t}\tbinom{m}{t}+1\ .
\end{split}
\end{equation*}
\end{itemize}
\end{proposition}

The bijections between the Farey sequence $\mathcal{F}_m$ and the
halfsequences of $\mathcal{F}\bigl(\mathbb{B}(2m),m\bigr)$,
presented in Theorem~\ref{prop:6}, allow us to describe the
properties of fractions from $\mathcal{F}_m$, analogous to those
of fractions from $\mathcal{F}\bigl(\mathbb{B}(2m),m\bigr)$,
presented in Proposition~\ref{prop:12}(ii,iii):

\begin{corollary} Let $m$ be an integer, $m>1$. The following combinatorial
identities hold for fractions from the standard Farey sequence
$\mathcal{F}_m$:
\begin{itemize}
\item[\rm(i)]
\begin{equation*}
\sum_{\substack{\frac{h}{k}\in\mathcal{F}_m:\\
\frac{0}{1}<\frac{h}{k}<\frac{1}{1}}}\ \sum_{1\leq s\leq
\left\lfloor\frac{m}{k}\right\rfloor}\tbinom{m}{s\cdot
h}\tbinom{m}{s\cdot k}=2^{2m-1}-2^m-\tfrac{1}{2}\tbinom{2m}{m}+1\
.
\end{equation*}
\item[\rm(ii)]
\begin{equation*}
\begin{split}
\quad&\phantom{=}\sum_{\substack{\frac{h}{k}\in\mathcal{F}_m:\\
\frac{0}{1}<\frac{h}{k}<\frac{1}{2}}}\ \sum_{1\leq s\leq
\left\lfloor\frac{m}{k}\right\rfloor}\tbinom{m}{s\cdot
k}\left(\tbinom{m}{s\cdot h}+\tbinom{m}{s\cdot(k-h)}\right)\\&=
\sum_{\substack{\frac{h}{k}\in\mathcal{F}_m:\\
\frac{1}{2}<\frac{h}{k}<\frac{1}{1}}}\ \sum_{1\leq s\leq
\left\lfloor\frac{m}{k}\right\rfloor}\tbinom{m}{s\cdot
k}\left(\tbinom{m}{s\cdot h}+\tbinom{m}{s\cdot(k-h)}\right)\\
&=2^{2m-1}-2^m-\tfrac{1}{2}\tbinom{2m}{m}- \sum_{1\leq t\leq
\left\lfloor\frac{m}{2}\right\rfloor}\tbinom{m}{2t}\tbinom{m}{t}+1\
.
\end{split}
\end{equation*}
\end{itemize}
\end{corollary}

\vskip 30pt

\section*{\normalsize References}

\noindent[1] D.~Acketa and J.~\v{Z}uni\'{c}, {\em On the number of linear partitions of the $(m,n)$-grid},
Inform. Process. Lett. {\bf 38} (1991), no.~3, 163--168.

\noindent[2] A.A.~Buchstab, {\em Teoriya Chisel} (in Russian) [{\em Number Theory}], Uchpedgiz, Moscow, 1960.

\noindent[3] S.B.~Gashkov and V.N.~Chubarikov, {\em Arifmetika. Algoritmy. Slozhnost' Vychisleniy, Third
edition} (in Russian) [Arithmetic. Algorithms. Complexity of Computation], Drofa, Moscow, 2005.

\noindent[4] R.L.~Graham, D.E.~Knuth and O.~Patashnik, {\em Concrete Mathematics. A Foundation for Computer
Science, Second edition}, Addison-Wesley, Reading Massachusetts, 1994.

\noindent[5] G.H.~Hardy and E.M.~Wright, {\em An Introduction to the Theory of Numbers, Fifth edition},
Clarendon Press, Oxford, 1979.

\noindent[6] A.O.~Matveev, {\em Pattern recognition on oriented matroids: Layers of tope committees}, {\tt
arXiv:math.CO/0612369}.

\noindent[7] A.O.~Matveev, {\em Relative blocking in posets}, J.~Comb.~Optim. {\bf 13} (2007), no.~4,
379--403.

\noindent[8] I.~Niven, H.S.~Zuckerman and H.L.~Montgomery, {\em An Introduction to the Theory of Numbers,
Fifth edition}, John Wiley \& Sons, Inc., New~York, 1991.

\noindent[9] J.J.~Tattersall, {\em Elementary Number Theory in Nine Chapters, Second edition}, Cambridge
University Press, Cambridge, 2005.

\end{document}